\documentclass{owrart}

\usepackage{latexsym, amsmath, amssymb, amsfonts, amscd}
\usepackage{amsthm}
\usepackage{t1enc}
\usepackage[mathscr]{eucal}
\usepackage{indentfirst}
\usepackage{fancyhdr}
\usepackage{fancybox}
\usepackage{enumerate}
\usepackage[all]{xy}

\theoremstyle{plain}
\newtheorem{thm}{Theorem}[section]

\newtheorem{prop}[thm]{Proposition}

\theoremstyle{definition}
\newtheorem{defi}[thm]{Definition}

\theoremstyle{remark}

\newtheorem{ep}[thm]{Example}

\newcommand{\lra}{\longrightarrow}

\newcommand{\rra}{\Rightarrow}

\newcommand{\Z}{\ensuremath{\mathbb Z}}
\newcommand{\C}{\ensuremath{\mathbb C}}
\newcommand{\R}{\ensuremath{\mathbb R}}

\newcommand{\cX}{\mathcal{X}}

\newcommand{\cG}{\mathcal{G}}
\newcommand{\cH}{\mathcal{H}}
\newcommand{\cS}{\mathcal{S}}

\newcommand{\ba}{\begin{eqnarray}}
   \newcommand{\na}{\end{eqnarray}}

\newcommand{\bt}{\mathbf{t}}                  
\newcommand{\bs}{\mathbf{s}}                  

\begin{document}

\begin{talk}[Henrique Bursztyn]{Chenchang Zhu}
{Morita equivalence of Poisson manifolds via stacky groupoids }
{Zhu, Chenchang}

\noindent

The aim of this talk is to present our program to define Morita equivalence in the
category of all Poisson manifolds via Morita equivalence of their
stacky symplectic groupoids. The talk is based on \cite{bz}. Early in \cite{xu}, Xu invented Morita
equivalence of Poisson manifolds with the inspiration from Rieffel's Morita
equivalence of $\C^*$-algebras. However it works only for
integrable Poisson manifolds, i.e. those Poisson manifolds who
process symplectic groupoids. A symplectic groupoid
\cite{wx} is a Lie groupoid $S\rra P$ with a symplectic
form $\omega$ on $S$ satisfying
\begin{equation}\label{form} pr_1^* \omega + pr_2^* \omega = m^* \omega, \end{equation}  on the set of composable
arrows $ S
\times_{P} S$ ($m$ is the multiplication on $S$). Then the base $P$ of the symplectic groupoid $S\rra P$ has
an induced Poisson structure such that the source map $\bs: S\to P$ is
a Poisson map and the target $\bt: S\to P$ is anti-Poisson. In fact there is a one-to-one correspondence between
integrable Poisson manifolds  and source-simply connected symplectic
groupoids.

Morita equivalence of Lie groupoids is well-studied and now
widely used in the theory of differentiable stacks. Roughly speaking, differentiable
stacks can be viewed as Lie groupoids up to Morita equivalence  (see for example
\cite{bx1}). Adding
compatible symplectic structure inside, \cite{xu} established Morita
equivalence of symplectic groupoids and proved further that 
 Poisson manifolds $P_1$ and
$P_2$ are Morita equivalent if and only if their symplectic groupoids
are Morita equivalent.  

Now \cite{tz} \cite{tz2} show that even a non-integrable Poisson
manifold processes a sort of symplectic groupoid $\cS \rra P$, but
$\cS$ is not anymore a manifold but an \'etale differentiable
stack\footnote{An \'etale differentiable stack is a differentiable
  stack presented by an \'etale Lie groupoid.  Careful readers find
  out that $\cS$ is presented by a groupoid and itself again is a
  groupoid over a manifold $P$. But these two groupoids are two
  different ones. In fact putting them together we have a Lie 2-groupoid
\cite{z:tgpd}.}
which processes a compatible symplectic form as in \eqref{form}.  Then the one-to-one
correspondence is extended to the set of all Poisson manifolds and
that of source-2-connected symplectic stacky groupoids (see Theorem \ref{main}). 

In our program, we first build Morita equivalence for stacky
groupoids, then we add compatible symplectic forms inside and build
Morita equivalence for symplectic stacky groupoids and hence for the base
Poisson manifolds.

\section{Stacky groupoids and their principal bundles}
We first say a few more words on the stacky groupoid $\cG
\rightrightarrows M$ we use. For an exact definition, we refer the reader to
\cite{z:tgpd}. The space of arrows $\cG$ is a differentiable stack,
and the space of objects $M$ is a manifold. It has $\bs$, $\bt$, $m$,
$e$, $i$ as
source, target, multiplication, identity, and inverse map
respectively, just as in the case of Lie groupoids. The only
difference now is that the multiplication is not strictly associative
but associative up to a 2-morphism $\alpha$ which satisfies a pentagon
condition. The same happens to all the other identities
we had before for Lie groupoids. Namely all these identities such as $(gh)k=g(hk)$,
$1g=g$, etc., do not hold strictly, but still hold up to something in a
controlled way. This `2'-phenomenon is new when we step into the world
of stacks. It will come back to haunt us all the time (for example Definition \ref{w-action}). The alternative
way is to work with Lie 2-groupoids which are essentially equivalent
to SLie groupoids \cite{z:tgpd}. We established
Morita equivalence of Lie 2-groupoids there. 

To shorten the notation, we call these stacky groupoids
{\em SLie groupoids}, and when $\cG$ is further an \'etale differentiable
stack,  a {\em W-groupoid}\footnote{The `W'
  comes from Alan Weinstein, who suggested this stacky approach to one of
  the authors.}. A {\em symplectic W-groupoid} is a
W-groupoid which has a compatible symplectic form as in \eqref{form}.

To build Morita equivalence, we first need the notion of principal
bundles of stacky groupoids.

\begin{defi}[SLie (W-)groupoid actions] \label{w-action}
Let $\cG$ be an SLie (W-)groupoid over $M$, $\cX$ differentiable
stack and $J: \cX\to M$ a smooth morphism. A right $\cG$-action
on $\cX$ is a smooth morphism
\[\Phi: \cX \times_{M} \cG \to \cX,\]
satisfies the following properties:
\begin{enumerate}
\item\label{itm:a-action} $\Phi \circ (\Phi \times id)=\Phi \circ (id
  \times m)$ holds up to a 2-morphism $a$;

\item $J \circ \Phi = \bs \circ pr_2$, where $pr_2:\cX \times_{M} \cG \to \cG$;

\item \label{itm:b-action} $\Phi \circ ( id \times ( e \circ J)) = id$
  holds up to a 2-morphism $b$.
\end{enumerate}
The 2-morphisms satisfy higher coherences, which roughly says that the
following diagrams commute:
\[
\xymatrix { & & ((x  g_1) g_2 ) g_3 \ar[dll] & & \\
(x g_1) (g_2 g_3) \ar[dr] & & & (x
(g_1 g_2)) g_3 \ar[ul] \\
& x  (g_1  (g_2 g_3)) \ar[r] & x ((g_1 g_2) g_3) \ar[ur] &}
\]
\[
\xymatrix{ & x(g\cdot 1) \ar[dl] \ar[dr] & \\
(xg)\cdot 1 \ar[rr] & & xg}
\]
\end{defi}

Given such an action, we can form a quotient stack $\cX/\cG$ as in
\cite{breen}. Unfortunately, the quotient stack is not always a
differentiable stack again. For this, we need principality of the action. 

Recall that an action $\Phi: X\times_M G \to X$ of a Lie groupoid  $G\rightrightarrows M$ on a
manifold $X$ is principal if and only if   $X/G$ is a manifold and
$pr_1 \times \Phi: X\times_M G \lra X\times_{X/G} X$ is an
isomorphism. We have the following definition:

\begin{defi}[{\bf Principal SLie (W-) groupoid bundles}] \label{prin}
Let $\cG \rra M$ be an SLie (-W) groupoid. A left
\textit{$\cG$-bundle over a differentiable stack $\cX$} is a
differentiable stack $\cX$ together with a smooth morphism
$\pi: \cX \to \cS$ and a right action $\Phi$ satisfying
\begin{equation}\label{proj-same}
\pi\circ \Phi = \pi\circ pr_2
\end{equation}
up to a 2-isomorphism $\alpha: \pi\circ pr_2 \to \pi \circ \Phi$. (Here
$pr_2:\cG\times_M\cX \to \cX$ is the natural projection.) The
2-isomorphism $\alpha$ satisfies a further coherence condition.

The bundle is \textit{principal} if  $\pi$ is a
surjective submersion  and
$$
pr_1 \times \Phi: \cX\times_{M} \cG \to \cX \times_{\cS} \cX
$$
is an isomorphism. Then the action $\Phi$ is also called \textit{principal}.
\end{defi}

\begin{ep} [A point as a principal $\Z$ bundle]
A point $pt$ is a principal $\Z$ bundle over the stack $B\Z$. The
action of $\Z$ on $pt$ is trivial, so it is not  principal  in the
classical sense. However, $pt$ is a principal $\Z$ bundle as in
Definition \ref{prin} because $pt \times_{B\Z} pt = \Z$ (see
\cite{bx1} for the definition of fibre product of differentiable stacks) and
\[ pt \times \Z \to pt \times_{B\Z} pt, \]
is an isomorphism of stacks.
\end{ep}

\begin{thm} \label{quotient} Let $\cG$ be an SLie (W-) groupoid. If $\pi: \cX \to \cS$ is a
$\cG$-principal bundle over $\cS$, then $\cX/\cG$ is a differentiable
stack and is isomorphic to
the base $\cS$. Moreover $\cX/\cG$ is presented by a Lie groupoid whose space of
arrows is $E_{\Phi}/G_1$ and whose space of objects is $ X_0$.
Here $X_1\rra X_0$ is a Lie groupoid presentation of $\cX$, $G_1\rra
G_0$ is that of $\cG$ and $E_\Phi$ is the H-S bibundle of the $\cG$-action $\Phi$.
\end{thm}

\section{Morita equivalence of SLie groupoids}

\begin{defi}[Morita equivalence of SLie groupoids]
Two SLie groupoids $\cG_1 \rightrightarrows M_1$ and $\cG_2
\rightrightarrows M_2$ are \textit{Morita equivalent} if there is
a differentiable stack $\cX$ and two smooth morphisms $J_i: \cX
\to \cG_i$ (moment maps) such that
\begin{enumerate}
\item  $J_1:\cX\to M_1$ is a right
principal $\cG_2$-bundle;
\item $J_2: \cX \to M_2$ is a left principal $\cG_1$-bundle;
\item $\Phi_2\circ (\Phi_1 \times id) = \Phi_1 \circ (id \times
  \Phi_2)$ holds up
to a 2-isomorphism $a$ which satisfies six higher coherence
conditions. 
\end{enumerate}
In this case we
call $\cX$ a $(\cG_1,\cG_2)$-{\em Morita bibundle}.
\end{defi}

It is simple to check that Morita equivalence is reflexive ($\cG$ itself is a
($\cG,\cG$)-Morita equivalence) and symmetric (use inverses to make right
actions into left and vice-versa). However transitivity is nontrivial
and we need to
use Theorem \ref{quotient}. 

Moreover we also have,

\begin{prop}
If two W-groupoids are Morita equivalent via Morita bibundle $\cX$,
then $\cX$ is an \'etale differentiable stack. 
\end{prop}

\begin{prop}
Two W-groupoids $\cG_i \rra M_i$ are Morita
equivalent via Morita bibundle $\cX$. If
$\cG_1\rra M_1$ is a Lie groupoid, then $\cX$ is a manifold and
$\cG_2\rra M_2$ is
also a Lie groupoid.
\end{prop}

Finally, two symplectic W-groupoids  $(\cG_1, \omega_1)
\rightrightarrows M_1$ and $(\cG_2, \omega_2)
\rightrightarrows M_2$ are \textit{Morita equivalent} if they are
Morita equivalent as SLie groupoids via a symplectic \'etale stack
$(\cX, \omega)$ satisfying
\[ pr_1^* \omega_1 + pr_2^* \omega = \Phi_1^* \omega, \quad
\text{on}\; \cG_1 \times_{M_1} \cX,  \]  where $\Phi_1$ is the action
of $\cG_1$ on $\cX$, and the same for $\omega$ and $\omega_2$.

\begin{thm}\label{main} \cite{tz2} For any symplectic W-groupoid $\cG\rightrightarrows M$, the base manifold $M$ has a unique Poisson structure such that the source map $\bs$ is Poisson. In this case, we call $\cG$ a  symplectic W-groupoid of the Poisson manifold $M$.

 On the other hand, for any Poisson manifold $M$, there are two
 symplectic groupoids $\cG(M)$ and $\cH(M)$ of $M$. $\cG(M)$ has
 2-connected source fibre and $\cH(M)$ has only 1-connected source fibre.
\end{thm}

\begin{defi}
Two Poisson manifolds  $M_1$ and $M_2$ are called {\em strongly Morita
equivalent} if $\cG(M_1)$  and $\cG(M_2)$ are Morita equivalent as
symplectic W-groupoids. Respectively, they are called {\em weakly
Morita equivalent} if $\cH(M_1)$ and $\cH(M_2)$ are Morita equivalent
as symplectic W-groupoids.
\end{defi} 

Strong Morita equivalence implies the weak one, and weak Morita equivalence coincides with the
classical one in \cite{xu} when applied to integrable Poisson
manifolds. But strong Morita equivalence is something new. For example, in \cite{xu}, with their usual symplectic forms,
$\R^2$ and the 2-sphere $S^2$ are Morita equivalent since all the
simply connected symplectic manifolds are Morita equivalent in the
classical sense. But they are {\em not} strongly Morita equivalent
because they have different $\pi_2$ groups. In
fact, only 2-connected symplectic manifolds are strongly Morita
equivalent to each other. We hope this $\pi_2$-phenomenon will help in
symplectic geometry, for example, in the aspect of preservation of prequantization.

\def\cprime{$'$} \def\cprime{$'$}

\end{talk}
\end{document}